\documentclass[a4paper,10pt]{amsart}
\usepackage{kubik}
\hypersetup{pdftitle={Local characterization of polyhedral spaces}%
,pdfauthor={Nina Lebedeva and Anton Petrunin}
}
\begin{document}
 \title{Local characterization of polyhedral spaces}
\author{Nina Lebedeva}
\address{N. Lebedeva\newline\vskip-4mm
Steklov Institute,
27 Fontanka, St. Petersburg, 
191023, Russia.
\newline\vskip-4mm
Math. Dept.
St. Petersburg State University,
Universitetsky pr., 28, 
Stary Peterhof, 
198504, Russia.}
\email{lebed@pdmi.ras.ru}
\author{Anton Petrunin}
\address{A. Petrunin\newline\vskip-4mm
Math. Dept. PSU,
University Park, PA 16802,
USA}
\email{petrunin@math.psu.edu}
\thanks{N.~Lebedeva was partially supported by RFBR grant 
14-01-00062.}
\thanks{A.~Petrunin was partially supported by NSF grant DMS 1309340.}

\date{}

\begin{abstract}
We show that  a compact length space is polyhedral if
a small spherical neighborhood of any point is conic. 
\end{abstract}
\maketitle

\section{Introduction}

In this note we characterize polyhedral spaces as the spaces where every point has a conic neighborhood.
Namely, we prove the following theorem;
see Section~\ref{sec:preliminaries} for all necessary definitions.

\begin{thm}{Theorem}\label{prop:poly-char}
A compact length space $X$ is polyhedral 
if and only if 
a neighborhood of each point $x\in X$ admits an open isometric embedding to Euclidean cone which sends $x$ to the tip of the cone.
\end{thm}

Note that we do not make any assumption on the dimension of the space.
If the dimension is finite
then the statement admits a simpler proof by induction; 
this proof is indicated in the last section.

A priori, it might be not clear why the space in the theorem
is even homeomorphic to a simplicial complex.
This becomes wrong if you remove word ``isometric'' from the formulation.
For example, there are closed 4-dimensional topological manifold 
which does not admit any triangulation, see \cite[1.6]{freedman}.

The Theorem~\ref{prop:poly-char} 
is applied in \cite{lebedeva},
where it is used to show that an Alexandrov space with the maximal number of extremal points is a quotient of $\RR^n$ by a cocompact properly discontinuous isometric action; see also \cite{lebedeva-0}.

\parbf{Idea of the proof.}
Let us cover $X$ by finite number of spherical conic neighborhood and consider its nerve, say $\mathcal{N}$.
Then we  map $\mathcal{N}$ barycentrically back to $X$.
If we could show that the image of this map cover whole $X$ 
that would nearly finish the proof. 
Unfortunately we did not manage to show this statement 
and have make a walk around; 
this is the only subtle point in the proof below.

\parbf{Acknowledgment.} 
We would like to thank
Arseniy Akopyan, 
Vitali Kapovitch,
Alexander Lytchak 
and Dmitri Panov
for their help.

\section{Definitions}\label{sec:preliminaries}

In this section we give the definition of polyhedral space of arbitrary dimension.
It seems that these spaces were first considered by Milka in \cite{milka};
our definitions are equivalent but shorter.

\parbf{Metric spaces.}
The distance between points $x$ and $y$ in a metric space $X$ will be denoted as $|x-y|$ or $|x-y|_X$.
Open $\eps$-ball centered at $x$ will be denoted as $B(x,\eps)$;
i.e.,
$$B(x,\eps)=\set{y\in X}{|x-y|<\eps}.$$
If $B=B(x,\eps)$ and $\lambda>0$ we
use notation $\lambda\cdot B$
as a shortcut for $B(x,\lambda\cdot\eps)$.

A metric space is called \emph{length space}
if the distance between any two points coincides with the infimum of lengths of curves connecting these points.

A \emph{minimizing geodesic} between points $x$ and $y$ will be denoted by $[xy]$.

\parbf{Polyhedral spaces.}
A length space is called \emph{polyhedral space}
if it admits a finite triangulation such that 
each simplex is (globally) isometric to 
a simplex in Euclidean space.

Note that according to our definition, 
the polyhedral space has to be compact.

\parbf{Cones and homotheties.}
Let $\Sigma$ be a metric space with diameter at most $\pi$.
Consider the topological cone
$K=[0,\infty)\times\Sigma/\sim$ where $(0,x)\sim(0,y)$ for every $x,y\in\Sigma$.
Let us equip $K$ with the metric defined by the rule of cosines;
i.e., for any $a,b\in[0,r)$ and $x,y\in\Sigma$ 
we have
$$|(a,x)-(b,y)|_K^2=a^2+b^2-2\cdot a\cdot b\cdot\cos |x-y|_\Sigma.$$
The obtained space $K$ 
will be called \emph{Euclidean cone over} $\Sigma$.
All the pairs of the type $(0,x)$ correspond to one point in $K$ 
which will be called the \emph{tip} of the cone.
A metric space which can be obtained in this way is called \emph{Euclidean cone}.

Equivalently, Euclidean cone can be defined as a metric space $X$ which admits a one parameter family of homotheties 
$m^\lambda\:X\to X$ for $\lambda\ge 0$ 
such that for any fixed $x, y\in X$
there are real numbers $\zeta$, $\eta$ and $\vartheta$
such that $\zeta,\vartheta\ge 0$, $\eta ^2\le \zeta\cdot\vartheta$ and
$$|m^\lambda(x)-m^\mu(y)|_X^2
=
\zeta\cdot \lambda^2+2\cdot \eta \cdot\lambda\cdot\mu+\vartheta\cdot \mu^2.$$
for any $\lambda,\mu\ge 0$.
The point $m^0(x)$ is the tip of the cone; 
it is the same point for any $x\in X$.

Once the family of homotheties is fixed,
we can abbreviate $\lambda\cdot x$ for $m^\lambda(x)$.

\parbf{Conic neighborhoods.}

\begin{thm}{Definition}
Let $X$ be a metric space, $x\in X$ and $U$ a neighborhood of $x$.
We say that $U$ is a \emph{conic} neighborhood
of $x$ if 
$U$ admits an open distance preserving embedding 
$\iota\:U\to K_x$ 
into Euclidean cone $K_x$
which sends $x$ to the tip of the cone.
\end{thm}

If $x$ has a conic neighborhood then the cone
 $K_x$ as in the definition will be called \emph{ the cone at $x$}.
Note that in this case $K_x$ is unique up to an isometry 
which sends the tip to the tip.
In particular, 
any conic neighborhood $U$ of $x$ admits an open distance preserving embedding $\iota_U\:U\to K_x$ which sends $x$ to the tip of $K_x$.
Moreover, 
it is easy to arrange that these embeddings commute with inclusions;
i.e., if $U$ and $V$ are two conic neighborhoods of $x$ and $U\supset V$
then the restriction of $\iota_U$ to $V$ coincides with $\iota_V$.
The later justifies that we omit index $U$ for the embedding $\iota\: U\to K_x$.

Assume $x\in X$ has a conic neighborhood and $K_x$ is the cone at $x$.
Given a geodesic $[xy]$ in $X$,
choose a point $\bar y\in \left[xy\right]\backslash\{x\}$ sufficiently close to $x$ and set
$$\log[xy]=\frac{|x-y|_X}{|x-\bar y|_X}\cdot\iota(\bar y)\in K_x.$$ 
Note that $\log[xy]$ does not depend on the choice of $\bar y$.

\section{Preliminary statements}

\begin{thm}{Definition}
Let $X$ be a metric space and $[px_1]$, $[px_2],\dots,[px_m]$ are geodesics in $X$.
We say that a neighborhood $U$ of $p$ 
\emph{splits} in the direction of the geodesics $[px_1]$, $[px_2],\dots,[px_k]$ 
if there is an open distance preserving map $\iota$ from $U$ to the product space $E\times K'$,
such that $E$ is a Euclidean space
and the inclusion
\[\iota(U\cap[px_i])\subset E\times \{o'\}\]
holds for a fixed $o'\in K'$ and any $i$.
\end{thm}

\parbf{Splittings 
and isometric copies of polyhedra.}
The following lemmas and the corollary are the key ingredients in the proof.

\begin{thm}{Lemma}\label{lem:splitting}
Let $X$ be a metric space, $p\in X$
and for each $i\zz\in\{1,\dots,k\}$
the ball
$B_i=B(x_i,r_i)$,  
forms a conic neighborhoods of $x_i$.
Assume $p\in B_i$
for each $i$.
Then  any conic neighborhood of $p$ splits in the direction of $[px_1],\dots,[px_k]$.
\end{thm}

In the proof we will use the following statement;
its proof is left to the reader.

\begin{thm}{Proposition}\label{prop:few-tips}
Assume $K$ is a metric space which admits cone structures with different tips
$x_1,\dots,x_k$.
Then $K$ is isometric to the product space $E\times K'$,
where $E$ is a Euclidean space 
and $K'$ is a cone with tip $o'$ and
$x_i\in  E\times \{o'\}$ for each $i$.
\end{thm}

\parit{Proof of Lemma~\ref{lem:splitting}.} 
Fix sufficiently small $\eps>0$.
For each point $x_i$,
consider point $x_i'\in[px_i]$ 
such that $|p-x_i'|=\eps\cdot|p-x_i|$.
Since $\eps$ is sufficiently small,
we can assume that $x_i'$ lies in the conic neighborhood of $p$.

Note that for the right choice of parameters close to 1,
the composition of homotheties with centers at $x_i$ and $p$ 
produce a homothety with center at $x_i'$
and these are defined in a fixed conic neighborhood of $p$.
(In particular it states that composition of homotheties of Euclidean space is a homothety;
the proof is the same. 
The parameters are assumed to be chosen in such a way that $x_i'$ stay fixed by the composition.) 

These homotheties can be extended to the cone $K_p$ at $p$
and taking their compositions we get the homotheties for all values of parameters with the centers at $\hat x'_i=\log[px_i']\in K_p$.
It remains to apply Proposition~\ref{prop:few-tips}.
\qeds

From the Lemma~\ref{lem:splitting}, we get the following corollary.

\begin{thm}{Corollary}\label{splitstructure}
Let $X$ be a compact length space and $x\in X$.
Suppose $B\zz=B(x,r)$ is a conic neighborhood of $x$ which
splits in the direction of $[px_1],\zz\dots,[px_k]$ and
$\iota\:B\hookrightarrow E \times K'$
be the corresponding embedding.
Then the image
$\iota(B)$ is a ball of radius $r$ centered at $\iota(x)\in E\times\{o'\}$.

In particular, for any point $q\in B$ such that
$|q-p|_X=\rho$ and
$\iota(q)\in E\times \{o'\}$
the ball $B(q,r-\rho)$ is a conic neighborhood of $q$.
\end{thm}

\begin{thm}{Lemma}\label{lem:polyhedron}
Let $B_i=B(x_i,r_i)$, $i\in\{0,\dots,k\}$
be balls in the metric space $X$.
Assume each $B_i$ forms a conic neighborhood of $x_i$
and $x_i\in B_j$ if $i\le j$.
Then $X$ contains a subset $Q$ which contains 
all $x_i$ and is isometric to a convex polyhedron.

Moreover the geodesics in $Q$ do not bifurcate in $X$;
i.e., if geodesic $\gamma\:[a,b]\to X$ lies in $Q$
and an other geodesics $\gamma'\:[a,b]\to X$ 
coincides with $\gamma$ on some interval then $\gamma'=\gamma$.
\end{thm}

To illustrate the second statement 
let us consider tripod $T$; 
i.e., 1-dimensional polyhedral space 
obtained from three intervals by gluing their left ends together.
Let $Q$ be the union of two segments in $T$.
Note that $Q$ forms a subset isometric to a real interval;
i.e., $Q$ is isometric to 1-dimensional convex polyhedron.
On the other hand, $Q$ does not satisfy the second condition since a geodesic can turn from $Q$ at the triple point.

\parit{Proof.}
To construct $Q=Q_k$ we apply induction on $k$ 
and use the cone structures on $B_i$ 
with the tip at $x_i$ 
consequently.

For the base case, $k=0$, we take $Q_0=\{x_{0}\}$.

By the induction hypothesis,
there is a set $Q_{k-1}$
containing all $x_0, \dots, x_{k-1}$.

Note that $B_k$ is strongly convex;
i.e., any minimizing geodesic with ends in $B_k$
lies completely in $B_k$.
In particular $Q_{k-1}\cap B_k$ is convex.
Since $x_i\in B_k$ for all $i<k$,
we may assume that $Q_{k-1}\subset B_k$.

Note that the homothety $m_k^\lambda$ with center $x_k$ and $\lambda\le 1$ is defined for all points in $B_k$.
Set 
\[Q_{k}=\set{m^\lambda_k(x)}{x\in Q_{k-1}\ \text{and}\  \lambda\le 1}.\] 
Since $Q_{k-1}$ is isometric to a convex polytope, 
so is $Q_k$.

To show that the geodesic $\gamma\:[a,b]\to X$ in $Q$
can not bifurcate, it is sufficient to show that 
if $a<c<b$ then 
a neighborhood of $p=\gamma(c)$
splits in the direction of $\gamma$.

The point $p$ can be obtained from from $x_0$
by a composition of homotheties
\[p=m_k^{\lambda_k}\circ\cdots\circ m_1^{\lambda_1}(x_0),\]
where
$0<\lambda_i\le1$.
Set $m=m_k^{\lambda_k}\circ\cdots\circ m_1^{\lambda_1}(x_0)$.
We can assume $r_0$ to be sufficiently small
so that $m$ is defined on
$B_0$.

By Lemma~\ref{lem:splitting}, 
$B_0$ splits in the directions of 
$[x_0x_1],\dots,[x_0x_k]$. 
Since $m$ rescales the distances by fixed factor,
a neighborhood of $p$ also splits.
Clearly the Euclidean factor in the image
$m(B_0)$
covers small neighborhood of $p$ in $Q$.
Since $\gamma$ runs $Q$,
a neighborhood of $p$
splits in the direction of $\gamma$.
\qeds

\section{The proof}

The proof of Theorem~\ref{prop:poly-char} is based on the following lemma.

\begin{thm}{Lemma}\label{lem:polytope=polytopes}
Assume a length space $X$ is covered by finite number of sets   
such that each finite intersection of these sets is isometric to a convex polytope. 
Then $X$ is a polyhedral space.
\end{thm}

\parit{Proof.}
It is sufficient to show that if any metric space $X$  (not necessary length metric space)
admits a cover as in the lemma then it admits a triangulation such that each simplex is isometric to a Euclidean simplex.

Denote by $V_1, \dots V_n$ the polytopes in the covering.
Let $m$ be the maximal dimension of $V_i$.

We will apply induction on $m$; the base case $m=0$ is trivial.

Now assume $m>0$.
Let $W_1\dots W_k$ denote all the faces of $V_1, \dots V_n$ of dimension at most $m-1$.
Note that the collection $W_1\dots W_k$ satisfies the assumption of the Lemma.
Therefore by induction hypothesis, $X'=\bigcup_i W_i$ admits the needed triangulation.

It remains to extend this triangulation to each of the $m$-dimensional polytopes which $X'$ cuts from $X$.
The later is generously left to the reader.
\qeds

\parit{Proof of Theorem~\ref{prop:poly-char}.}
We need to show the ``if'' part;
the ``only if'' part is trivial.

Fix a finite cover of $X$ by open balls $B_i=B(x_i,r_i)$,
 $i\in\{0,\dots,n\}$
such that for each $i$, 
the ball $5\cdot B_i$
forms a conic neighborhood of $x_i$.

Given $i\in\{0,\dots,n\}$ and $z\in X$ set
$$f_i(z)=|x_i-z|_X^2-r_i^2.$$ 
Clearly $f_i(z)<0$ if and only if  $z\in B_i$. 

Set 
$$f(z)=\min_i \{f_i(z)\}.$$
It follows that $f(z)<0$ for any $z\in X$. 

Consider \emph{Voronoi domains} $V_i$ for the functions $f_i$; 
i.e.,
$$V_i=\set{z\in X}{f_i(z)\le f_j(z)\ \text{for all}\ j}.$$
From above we get that $V_i\subset B_i$ for each $i$.%
\footnote{It also follows that $V_i$ forms a \emph{strongly convex subset} of $X$; 
i.e., any minimizing geodesic in $X$
with ends in $V_i$ lies completely in $V_i$.
This property is not needed in our proof, 
but it is used in the alternative proof;
see the last section.}

Given a subset $\sigma\subset\{0,\dots,n\}$
set 
\[V_\sigma=\bigcap_{i\in\sigma} V_i.\]
Note that $V_{\{i\}}=V_i$ for any $i\in\{0,\dots,n\}$.

Let $\mathcal N$ be the \emph{nerve} of the covering $\{V_i\}$;
i.e., $\mathcal N$ is the abstract simplicial complex 
with $\{0,\dots,n\}$ as the set of vertexes
and such that a subset $\sigma\zz\subset\{0,\dots,n\}$ 
forms a simplex in $\mathcal N$ 
if and only if
$V_\sigma\ne\emptyset$.

Let us fix a simplex $\sigma$ in $\mathcal N$.
While $\sigma$ is fixed, we may assume without loss of generality
that $\sigma=\{0,\dots,k\}$ for some $k\le n$ 
and $r_0\le r_1\le\dots r_k$.
In particular $2\cdot B_i\ni x_0$ for each $i\le k$.

From above $V_\sigma\subset B_0$.
Since $5\cdot B_i$ is a conic neighborhood of $x_i$
and $2\cdot B_i\ni x_0$ for each $i\in \sigma$,
we can apply Lemma~\ref{lem:splitting} for the balls 
$5\cdot B_0,\dots,5\cdot B_k$.
Denote by $h\:5\cdot B_0\hookrightarrow E\times K$
the distance preserving embedding provided by this lemma.
We can assume that the Euclidean factor $E$
has minimal possible dimension; i.e., the images $h(B_0\cap [x_0x_i])$ span whole $E$.
In this case the projection of $h(V_\sigma)$ on $E$ is a one-point set, say $\{z\}$.
Denote by $x_\sigma\in B_0$ the point such that $h(x_\sigma)=z$.
Set $r_\sigma=r_0$ and $B_\sigma=B(x_\sigma,r_\sigma)$.
(The point $x_\sigma$ plays the role of \emph{radical center} of the collection of balls $\{B_i\}_{i\in\sigma}$.)

According to Corollary~\ref{splitstructure} 
the ball $4\cdot B_\sigma$ 
forms a conic neighborhood of $x_\sigma$.
Clearly $B_\sigma\supset V_\sigma$.

Let $\phi$ and $\psi$ be faces of $\sigma$;
in other words, $\phi$ and $\psi$ are subsets in $\sigma=\{0,\dots,k\}$.
Set $i=\min\phi$ and $j=\min\psi$.
Assume $i\ge j$, in this case $r_\phi=r_i\ge r_j =r_\psi$.
From above we get $x_\phi\in B_i$, $x_\psi\in B_j$ and $x_j\in 2\cdot B_i$. 
Therefore $x_\psi\in 4\cdot B_\phi$.

Therefore Lemma~\ref{lem:polyhedron} provides a subset, say $Q_\sigma$ isometric to a convex polyhedron and contains all
$x_\phi$ for $\phi\subset\sigma$.

It remains to show 
\begin{enumerate}[(a)]
\item\label{a} $X=\bigcup_{\sigma}Q_\sigma$,  where the union is taken for all the simplices $\sigma$ in $\mathcal N$.
\item\label{b} The intersection of arbitrary collection of $Q_\sigma$ is isometric to a convex polytope.
\end{enumerate}
Once (\ref{a}) and (\ref{b}) are proved,
Lemma \ref{lem:polytope=polytopes} will finish the proof.

Part (\ref{b}) follows since the geodesics in $Q_\sigma$
do not bifurcate in $X$; see Lemma~\ref{lem:polyhedron}.

Given $p\in X$,
set
$$
\sigma(p)
=
\set{i\in\{0,\dots,n\}}{p\in V_i}.
$$
Note that $\sigma(p)$ forms a simplex in $\mathcal{N}$
and $p\in V_{\sigma(p)}$.
Therefore $p\in B_{\sigma(p)}$.

Recall that $B_{\sigma(p)}$ forms a conic neighborhood of $x_{\sigma(p)}$. 
If $p\ne x_{\sigma(p)}$ then moving $p$ away from $x_{\sigma(p)}$ in the radial direction keeps the point in  $V_{\sigma(p)}$
till the moment it hits a new Voronoi domain, say $V_j$ with $j\notin \sigma(p)$.
Denote this end point by $p'$.
In other words, $p'$ is the point such that 
\begin{enumerate}[(i)]
\item $p$ lies on the geodesic $[x_{\sigma(p)} p']$;
\item $p'\in V_i$ for any $i\in \sigma(p)$;
\item the distance $|x_{\sigma(p)} - p'|_X$ takes the maximal possible value.
\end{enumerate}

{
\begin{wrapfigure}{r}{38mm}
\begin{lpic}[t(-3mm),b(4mm),r(0mm),l(0mm)]{pics/p123(1.0)}
\lbl[t]{2,0;$x_2$}
\lbl[b]{36,26.5;$x_1$}
\lbl[b]{2,26.5;$x_3$}
\lbl[tr]{25,20.3;$p_0$}
\lbl[tr]{14,19;$p_1$}
\lbl[bl]{10,29;$p_2=x_{\{1,2,3\}}$}
\lbl[tl]{20.5,14;$x_{\{1,2\}}$}
\lbl{32,7;$V_1$}
\end{lpic}
\end{wrapfigure}

Start with arbitrary point $p$
and consider the recursively defined sequence $p=p_0,p_1,\dots$
such that $p_{i+1}=p_i'$.

Note that $\sigma(p)$ forms a proper subset of $\sigma(p')$.
It follows that the sequence $(p_i)$ terminates after at most $n$ steps;
in other words $p_k=x_{\sigma(p_k)}$ for some $k$.

In particular $p_k\in Q_{\sigma(p_k)}$.
By construction it follows that $p_i\in Q_{\sigma(p_k)}$
for each $i\le k$.
Hence $p\in Q_{\sigma(p_k)}$; i.e., (\ref{a}) follows.
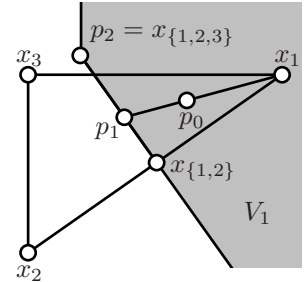
\qeds
}

\section{Final remarks}\label{sec:remarks}

\parbf{Finite dimensional case.}
Let $X$ be a compact length space 
such that each point $x\in X$
admits a conic neighborhood.

Note that from Theorem~\ref{prop:poly-char},
it follows in particular that dimension of $X$ is finite.
If we know a priori the dimension (topological or Hausdorff) of $X$ is finite
then one can build an easier proof using induction on the dimension 
which we are about to indicate.

Consider the Voronoi domains $V_i$ as in the beginning of proof of Theorem~\ref{prop:poly-char}.
Note that all $V_i$ are convex and 
\[\dim V_{\{i,j\}}<\dim X\] if $i\ne j$.

By induction hypothesis we can assume that all $V_{\{i,j\}}$ are polyhedral spaces. 
Cover each $V_{\{i,j\}}$ by 
isometric copies of convex polyhedra satisfying Lemma~\ref{lem:polytope=polytopes}.
Applying the cone construction with center $x_i$ over these copies in $V_{\{i,j\}}$ for all $i\ne j$,
we get a covering of $X$ by a finite number of 
copies of convex polyhedra 
such that all their finite intersections are isometric to convex polyhedra.
It remains to apply Lemma~\ref{lem:polytope=polytopes}.

\parbf{Spherical and hyperbolic polyhedral spaces.}
Analogous characterization holds 
for spherical and hyperbolic polyhedral spaces. 
One needs to use spherical 
and hyperbolic rules of cosine in the definition of cone;
after that proof goes without any changes.

\parbf{Locally compact case.}
One may define polyhedral space 
as a complete length space which admits a locally finite triangulation 
such that each simplex is isometric to 
a simplex in Euclidean space.

In this case a locally compact length space is polyhedral
if every point admits a conic neighborhood.
The proof is the same.

\parbf{One more curvature free result.}
Our result is curvature free ---
we do not make any assumption on the curvature of $X$.
Besides our theorem,
we are aware about only one statement of that type  
--- the polyhedral analog of Nash--Kuiper theorem.
It states that any distance nonexpanding map from $m$-dimensional polyhedral space to the Euclidean $m$-space
can be approximated by a piecewise distance preserving map to the Euclidean $m$-space.
In full generality this result was proved recently by Akopyan \cite{akopyan}, 
his proof is based on earlier results obtained
by Zalgaller \cite{zalgaller} and Krat \cite{krat}.
Akopyan's proof is sketched in the lecture notes of the second author \cite{PY}.

\end{document}